\documentclass{amsart}

\usepackage{graphicx}
\usepackage[a4paper]{geometry}
\usepackage{booktabs}

\newcommand{\pHmac}{\ensuremath{\text{pH}_{\textrm{macro}}}}
\newcommand{\R}{\ensuremath{\mathbb{R}}}
\newcommand{\dx}[1]{\,\mathrm{d}#1}

\begin{document}

\title[Computations of corrosion front in sewer pipes]{Macroscopic corrosion
front computations of sulfate attack in sewer pipes based on a micro-macro
reaction-diffusion model}

%    author one information
\author{Vladim\'ir Chalupeck\'y}
\address{Institute of Mathematics for Industry, Kyushu University, Japan}
\curraddr{}
\email{chalupecky@imi.kyushu-u.ac.jp}
\thanks{}

%    author two information
\author{Tasnim Fatima}
\address{Centre for Analysis, Scientific Computing, and Applications, Department
of Mathematics and Computer Science, Eindhoven University of Technology, The
Netherlands}
\curraddr{}
\email{t.fatima@tue.nl}
\thanks{}

%    author three information
\author{Jens Kruschwitz}
\address{Strassen.NRW, Gelsenkirchen, Germany}
\curraddr{}
\email{jens.kruschwitz@strassen.nrw.de}
\thanks{}

%    author four information
\author{Adrian Muntean}
\address{Centre for Analysis, Scientific Computing, and Applications, Department
of Mathematics and Computer Science, Institute for Complex Molecular Systems,
Eindhoven University of Technology, The Netherlands}
\curraddr{}
\email{a.muntean@tue.nl}
\thanks{}


\date{}

\dedicatory{}

\begin{abstract}
  We consider a two-scale reaction diffusion system able to capture the
  corrosion of concrete with sulfates. Our  aim here is to define and compute
  two macroscopic corrosion indicators: typical pH drop and gypsum profiles.
  Mathematically, the system is coupled, endowed with micro-macro transmission
  conditions, and posed on two different spatially-separated scales: one
  microscopic  (pore scale) and  one macroscopic (sewer pipe scale). We use a
  logarithmic expression to compute values of pH from the volume averaged
  concentration of sulfuric acid which is obtained by resolving numerically the
  two-scale system (microscopic equations with direct feedback with the
  macroscopic diffusion of one of the reactants). Furthermore, we also evaluate
  the content of the main sulfatation reaction (corrosion) product---the
  gypsum---and point out numerically a persistent kink in gypsum's concentration
  profile. Finally, we illustrate numerically the position of the free boundary
  separating corroded from not-yet-corroded regions. 
\end{abstract}

\maketitle

\section{Introduction}\label{sec:intro}

\subsection{Background on sulfate corrosion}

Often in service-life predictions of concrete structures (e.g., sewer systems),
the effects of chemical and biological corrosion processes are fairly neglected.
In sewer systems and wastewater treatment facilities, where high concentrations
of hydrogen sulfide, moisture, and oxygen are present in the atmosphere, the
deterioration of concrete is caused mainly by biogenic acids. The so-called
microbially-induced concrete corrosion in sewer systems has been a serious
unsolved problem\footnote{There are a lot of financial implications if you want
to change the network of pipes in a city like Fukuoka. Our statement here is
that questions like {\em Why changing the pipes if corrosion is not so strong
yet and therefore the mechanics structure of the network can/could still hold
for 5 more years?} can be addressed in a rigorous mathematical multiscale
framework. Such an approach would allow a good understanding and prediction at
least of extreme situations.} for long time. The presence of microorganisms such
as fungi, algae or bacteria can induce formation of aggressive biofilms on
concrete surfaces. Particularly, the sulfuric acid that causes corrosion of
sewer crowns is generated by such a complex microbial ecosystem especially in
hot environments. The precise role of microorganisms in the context of sulfates
attack on concrete (here we focus on sewer pipes) is quite complex and is
therefore less understood from both experimental and theoretical points of view;
see, e.g., the experimental studies \cite{BielefeldtEtAl2010} (optimum pH and
growth kinetics of four relevant bacterial strains), \cite{IslanderEtAl1991}
(characteristics of the crown microbial system), \cite{Masch_Bock1994}
(microbiologically influenced corrosion of natural sandstone), \cite{Okabe}
(succession of sulfur-oxidizing bacteria in the bacterial community on corroding
concrete), \cite{Parker1945} (isolation of {\em Thiobacillus thiooxidans}),
\cite{S_B1984} (Hamburg sewers),\cite{Yongsiri2004} (air-water transfer of
hydrogen sulfide). As a consequence of this, an accurate large-time forecast of
the penetration of the sulfate corrosion front is very difficult to obtain. 

We want to stress the fact that concrete, in spite of its strong heterogeneity,
is mechanically a well-understood material with known composition. Also, the
cement (paste) chemistry is well understood. However, all cement-based materials
(including concrete) involve a combination of ``heterogeneous multi-phase
material'', ``multiscale chemistry'',  ``multiscale transport'' (flow,
diffusion, ionic fluxes, etc.), and ``multiscale'' mechanics. Having in view
this complexity, such materials are \emph{sensu stricto} very difficult to
describe, to analyze mathematically, and last but not least, to deal with
numerically. We expect that only after the multiscale aspects of such materials
are handled properly, good predictions of the large-time behavior may be
obtained. This is our path to addressing this corrosion scenario that is often
referred to as \emph{the sulfatation problem}.

Before closing these background notes, let us add some remarks \cite{Raats} on a
closely-related topic of acid sulphate soils\footnote{Compared to concrete,
soils are much easier to handle. Their mechanics is simpler and their chemistry
is often rudimentary, if any.}, which might attract the attention of the
multiscale research in the near future. Acid sulphate soils are an important
class of soils worldwide. Particularly in coastal areas, sediments often contain
a large amount of iron sulfide ($\mathrm{FeS}$ and/or $\mathrm{FeS_2}$). When by
drainage the sediment is exposed to air, this iron sulfide will oxidize to iron
sulphate ($\mathrm{FeSO_4}$). As long as the sediment still contains calcium
carbonate ($\mathrm{CaCO_3}$), the $\mathrm{FeS}$ will react with it, resulting
in gypsum ($\mathrm{CaSO_4}$) and iron oxide ($\mathrm{Fe_2O_3}$). Gypsum, being
much more soluble than calcium carbonate, will tend to leach to the ground- and
surface-waters. If the $\mathrm{Fe_SO_4}$ is no longer removed by reactions with
$\mathrm{CaCO_3}$ or other materials, it will tend to accumulate, resulting in a
drop of the pH below 4. The problematic acid sulphate soil then will have become
a reality. Acid sulphate soils were first described in 1886 by the Dutch Chemist
Jacob Maarten van Bemmelen, in connection with problems arising in the
Haarlemmermeer Polder. Much of the work in this direction is/was done in the
tropics, including Indonesia, Vietnam, and Australia; see, e.g., \cite{CSIRO}.

\subsection{Objectives and structure of the paper}

In order to be able to tackle the biophysics of the problem at a later stage,
eventually coupled with the mechanics of the concrete and the actual capturing
of the macroscopic fracture initiation, we focus here on a much simpler setting
modeling the multiscale transport and reaction of the active chemical species
involved in the sulfatation process. Therefore, the approach and results
reported here are only preliminary. 

Our main objective is twofold: using a multiscale reaction-diffusion system for
concrete corrosion (that allows for feedback between micro and macro scales), 
\begin{itemize}
  \item calculate pH profiles and detect the eventual presence of ``sudden'' pH
    drops; 
  \item extract from gypsum concentration profiles the approximate position of macroscopic
    corrosion fronts.
\end{itemize}

% \tableofcontents

In Section~\ref{sec:chem}, we present the reaction mechanisms taking place in
sewer pipes. In Section~\ref{sec:prob}, we give a mathematical description of
the problem and we set a two-scale PDE-ODE system. We briefly comment on a few
mathematical properties of the model. In Section~\ref{sec:num}, we approximate a
macroscopic pH numerically using a multiscale FD scheme and comment briefly on
the numerical results.

\section{A few notes on the involved chemistry}\label{sec:chem}

Our model includes two important features:
\begin{itemize}
  \item continuous transfer of $\mathrm{H_2S}$ from water to air phase and
    vice versa;
  \item fast production of \emph{gypsum} at solid-water interface.
\end{itemize}
We incorporate the Henry's law to model the transfer of $\mathrm{H_2S}$ from the
water to the air phase and vice versa \cite{Balls,Yongsiri2004}. The production
of gypsum at the solid-water interface is modeled by a non-linear reaction rate,
given by \eqref{reaction-rate}.

\begin{figure}[ht]
  \centering
  \includegraphics[width=\textwidth]{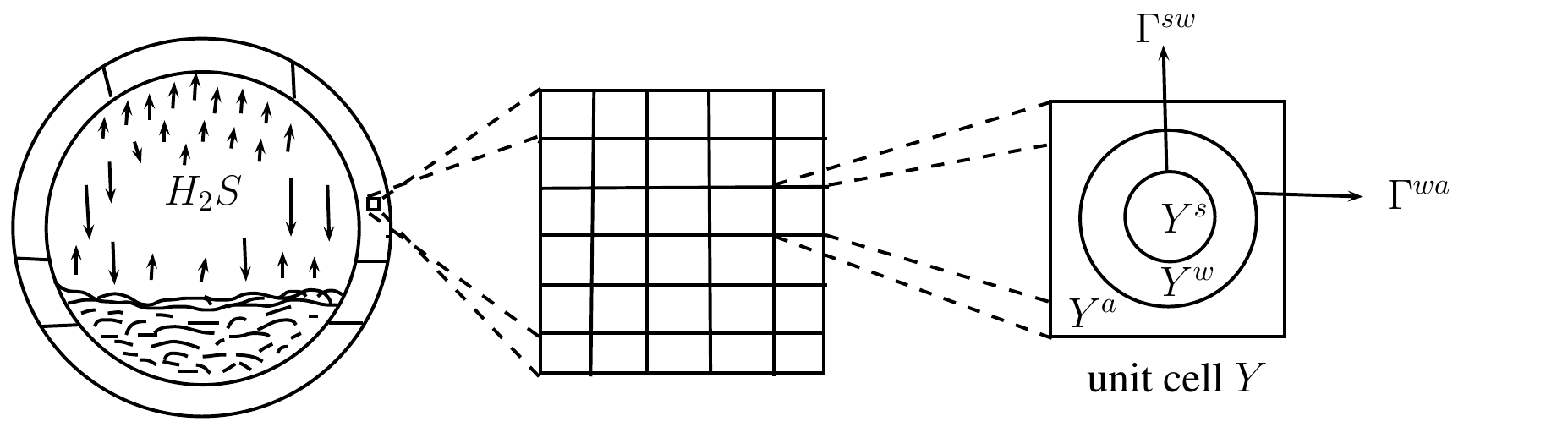}
  \caption{Left: Cross-section of a sewer pipe. Middle: Mesoscopic periodic
  approximation of a REV.  Right: Our concept of pore geometry
  (microstructure).}
  \label{fig:1}
\end{figure}

There are many variants of severe attack to concrete in sewer pipes which
influence the performance of concrete structure depending on the intensity of
the reactions, the environment, and the turbulence of the wastewater
\cite{Yong1}. We focus here on the most aggressive one, namely we consider the
following reaction mechanisms causing sulfatation, viz.
\begin{align}
  \label{eq:reaction-equations-1}
  \mathrm{10\,H^+} + \mathrm{SO_4^{2-}} + \text{org. matter} & \rightarrow
  \mathrm{H_2S(aq)} + \mathrm{4\,H_2O} + \text{oxidized matter}\\
  \label{eq:reaction-equations-2}
  \mathrm{H_2S(aq)} + \mathrm{2\,O_2} &\rightarrow \mathrm{2H^+} +
  \mathrm{SO_4^{2-}}\\
  \label{eq:reaction-equations-3}
  \mathrm{H_2S(aq)} &\rightleftharpoons \mathrm{H_2S(g)}\\
  \label{eq:sulfatation}
  \mathrm{2\,H_2O} + \mathrm{H^+} + \mathrm{SO_4^{2-}} + \mathrm{CaCO_3}
  &\rightarrow \mathrm{CaSO_4 \cdot 2\,H_2O} + \mathrm{HCO_{3-}}
\end{align}
Reaction \eqref{eq:reaction-equations-3} is typically a surface reaction taking
place as soon as water and air phases meet together. It plays an important role
in transferring the $\mathrm{H_2S}$ from the air phase to the liquid phase where
the corrosion actually takes place. For modeling details such as a Henry-like
``reaction'' mechanism, we refer the reader to \cite{tasnim1,mbp} and references
cited therein.

\section{Multiscale description of the sulfatation problem}\label{sec:prob}

We assume that the geometry of our concrete sample (porous medium) consists of a
system of pores periodically distributed inside a three-dimensional cube
$\Omega:=[a,b]^3$ with $a,b\in\R$ and $b>a$. The exterior boundary of $\Omega$
consists of two disjoint, sufficiently smooth parts: the Neumann boundary
$\Gamma^N$ and the Dirichlet boundary $\Gamma^D$. We assume that the pores in
concrete are made of stationary water film, air and solid parts in different
ratios depending on the local porosity. The reference pore, say $Y:=[0,1]^3$,
has three pair-wise disjoint domains $Y^s$, $Y^w$ and $Y^a$ with smooth
boundaries $\Gamma^{sw}$ and $\Gamma^{wa}$ as shown in Fig.~\ref{fig:1} such
that
\[
    Y = \bar{Y}^s \cup \bar{Y}^w \cup \bar{Y}^a.
\]
We refer the reader to \cite{tasnim1,ChalupeckyEtAl2010} for more description of
the multiscale geometry of the porous material. For a single scale (macroscopic)
approach of a sulfatation scenario, we refer the reader to \cite{natalini}, e.g. 

We consider a two-scale system of PDEs and one ODE for unknown functions
$u_1:\Omega\times (0,T)\rightarrow \R$, $u_k:\Omega\times Y^w\times (0,T)
\rightarrow \R$, $k\in\{2,3\}$, and $u_4:\Omega\times\Gamma^{sw}\times
(0,T)\rightarrow \R$ where $(0,T)$ is the time interval. The model under
consideration is derived by formal homogenization using different scalings of
the diffusion coefficients in \cite{tasnim1} (see also
\cite{MunteanChalupecky2011}) and is given by
\begin{subequations}\label{eq:two-scale-system}
\begin{align}
  \label{eq:two-scale-system-u1}
    \partial_t u_1 - d_1 \Delta u_1 &= - B\left(H u_1 - \int_{\Gamma^{wa}} u_2
    \dx{\gamma_y} \right),& \text{ in } &\Omega\times (0,T),\\
  \label{eq:two-scale-system-u2}
    \beta_2 \partial_t u_2 - \beta_2 d_2 \Delta_y u_2 &= -\Phi_2 k_2 u_2 + \Phi_3 k_3 u_3,& \text{ in } &\Omega\times Y^w \times (0,T),\\
  \label{eq:two-scale-system-u3}
    \beta_3 \partial_t u_3 - \beta_3 d_3 \Delta_y u_3 &=  \Phi_2 k_2 u_2 - \Phi_3 k_3 u_3,& \text{ in } &\Omega\times Y^w \times (0,T),\\
  % \label{eq:two-scale-system-w4}
    % \beta_4 \partial_t w_4 - \beta_4 d_4 \Delta_y w_4 &= k_1 w_1,
    % \mbox{ in } \Omega\times Y^w \times (0,T)\\
  \label{eq:two-scale-system-u4}
    \beta_4 \partial_t u_4 &= \Phi_4 k_4 \eta(u_3,u_4),& \text{ in } &\Omega\times \Gamma^{sw}\times (0,T),
\end{align}
\end{subequations}
where $u_1$ denotes the concentration for $\mathrm{H_2S}$ gaseous species, $u_2$
for $\mathrm{H_2S}$ aqueous species, $u_3$ for $\mathrm{H_2SO_4}$, and $u_4$ for
\emph{gypsum} at $\Gamma^{sw}$. The water film is taken here to be stationary. A
detailed modeling of the role of water is still open, see, e.g.,
\cite{Raats2,beddoe1,tixier}. $\Delta$ without subscript denotes the Laplace
operator with respect to macroscopic variable $x$ and $\Delta_y$ with respect to
microscopic variable $y$. $\dx{\gamma_y}$ represents the differential over the
surface $\Gamma^{wa}$. $\beta_k>0$, $k\in\{2,3,4\}$, represents the ratio of the
maximum concentration of the $k$-th species to the maximum concentration of
$\mathrm{H_2SO_4}$, $d_i>0$, $i\in\{1,2,3\}$, are the diffusion coefficients,
$B$ is a dimensionless Biot number which gives the mass transfer rate between
water and air phases, and $k_j:Y\rightarrow R$, $j\in\{2,3,4\}$, are functions
modeling the reaction rate ``constants''. $\Phi_k$ ($k\in\{2,3,4\}$) are
Damk\"ohler numbers corresponding to three distinct chemical mechanisms
(reactions). They are dimensionless numbers comparing the characteristic time of
the fastest transport mechanism (here, the diffusion of $\mathrm{H_2S}$ in the
gas phase) to the characteristic timescale of the $k$-th chemical reaction.

The system \eqref{eq:two-scale-system} is supplemented with initial and boundary
conditions, which read as
\begin{align}
  \label{eq:two-scale-system-ic1}
  u_1(x,0)   &= u^0_1(x),& \text{ on } &\Omega \times (0,T),\\
  \label{eq:two-scale-system-ic2}
  u_k(x,y,0) &= u^0_k(x,y),& \text{ on } &\Omega \times Y^w \times (0,T),\ k\in\{2,3\},\\
  \label{eq:two-scale-system-ic3}
  u_4(x,y,0) &= u^0_4(x,y),& \text{ on } &\Omega \times \Gamma^{sw} \times (0,T),\\
  \label{eq:two-scale-system-bc1}
  u_1 &= u_1^D,& \text{ on } &\Gamma^D \times (0,T),\\
  \label{eq:two-scale-system-bc2}
  n_N\cdot(d_1\nabla u_1) &= 0,& \text{ on } &\Gamma^N \times (0,T),\\
  \label{eq:two-scale-system-bc3}
  n_{wa}\cdot(d_2\nabla_y u_2 ) &= B \left(H u_1 - \int_{\Gamma^{wa}} u_2
  \dx{\gamma_y}\right),& \text{ on } &\Omega \times \Gamma^{wa} \times (0,T),\\
  \label{eq:two-scale-system-bc4}
  n_{sw}\cdot(d_2\nabla_y u_2 ) &= 0,& \text{ on } &\Omega \times \Gamma^{sw} \times (0,T),\\
  \label{eq:two-scale-system-bc5}
  n_{wa}\cdot(d_3\nabla_y u_3) &= 0,& \text{ on } &\Omega\times \Gamma^{wa} \times (0,T),\\
  \label{eq:two-scale-system-bc6}
  n_{sw}\cdot(d_3\nabla_y u_3) &= -\Phi_3 \eta(u_3, u_4),& \text{ on } &\Omega
  \times \Gamma^{sw} \times (0,T),
  % \label{eq:two-scale-system-bc7}
  % n\cdot(-d_4\partial_y w_4 )&= 0, \qquad \mbox{ on } \Omega\times \Gamma^{sw}\times (0,T),\\
  % \label{eq:two-scale-system-bc8}
  % n\cdot( -d_4\partial_y w_4 )&= 0, \qquad \mbox{ on } \Omega\times\Gamma^{wa}\times (0,T).
\end{align}
where $n_N$ denotes the outward unit normal vector to $\partial\Omega$ along
$\Gamma^N$, and $n_{wa}$ and $n_{sw}$ denote the outward unit normal vectors to
$Y^w$ along $\Gamma^{wa}$ and $\Gamma^{sw}$, respectively. Note that the
``information'' at the micro-scale is connected to the macro-scale situation via
the right-hand side of \eqref{eq:two-scale-system-u1} and via the
\emph{micro-macro boundary condition} \eqref{eq:two-scale-system-bc3}. It is
also worth noticing that all involved parameters (except for $H$, $d_3$ and $B$)
contain microscopic information. The coefficients $d_3$ and $B$ are effective
ones (see \cite{tasnim1,eck} for a way of calculating them), while $H$ can be
read off from existing macroscopic experimental data.

We consider the following form of the reaction rate $\eta$ at the interface
$\Gamma^{sw}$
\begin{equation}\label{reaction-rate}
  \eta(\alpha,\beta)=
  \begin{cases}
    \alpha^p(\bar{\beta}-\beta)^q, & \text{ if } \alpha \geq 0, \beta \geq 0,\\
    0, & \mbox{otherwise},
  \end{cases}
\end{equation}
where $\bar{\beta}$ is a known maximum concentration of gypsum at $\Gamma^{sw}$
and $p\geq 1,q\geq1$ are partial orders of reaction. For more modeling
possibilities of choosing $\eta$, see \cite{tasnim3}. It is worth noting that
production terms like
\[
B\left( H u_1(t,x) -  \int_{\Gamma^{wa}}u_2(t,x,y) \dx{\gamma_y}\right)
\]
are usually referred in the literature as Henry's or Raoult's law, where $H>0$
is known Henry's constant.

We refer the reader to \cite[Theorem~3]{ChalupeckyEtAl2010} for statements
regarding the global existence and uniqueness a weak solution to problem
\eqref{eq:two-scale-system}--\eqref{eq:two-scale-system-bc6} (see also
\cite{MunteanNeuss2010} for the analysis on a closely related problem).

\section{Simulation at a macroscopic pH scale. Capturing free boundaries}\label{sec:num}

In this section the model \eqref{eq:two-scale-system} is applied to the
simulation of the acid corrosion due to a microbiotical layer on a cement
specimen. We focus on extracting the position of the corrosion front and on the
acid reaction, which we use to obtain macro-scale profiles of pH. Both of these
results can be compared to experimental data published, e.g., in
\cite{IslanderEtAl1991,Okabe}.

\begin{figure}[ht]
  \begin{center}
    \includegraphics[width=0.48\textwidth]{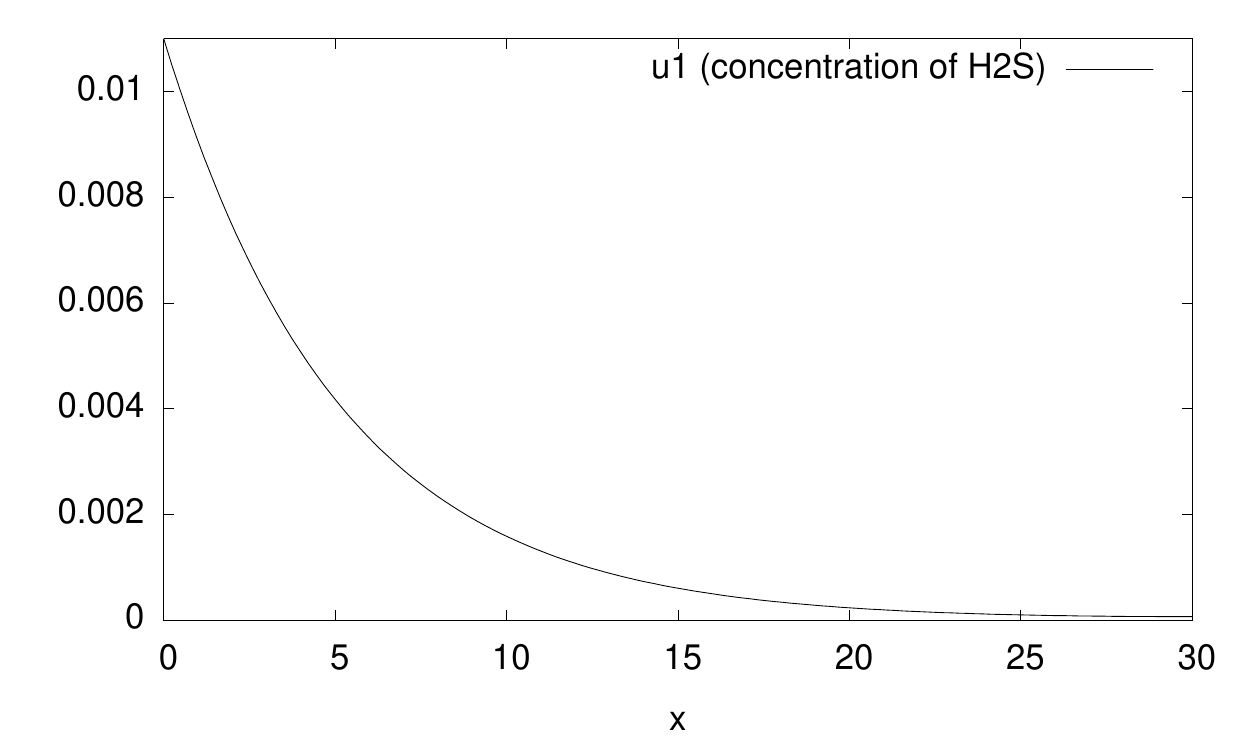}
    \includegraphics[width=0.48\textwidth]{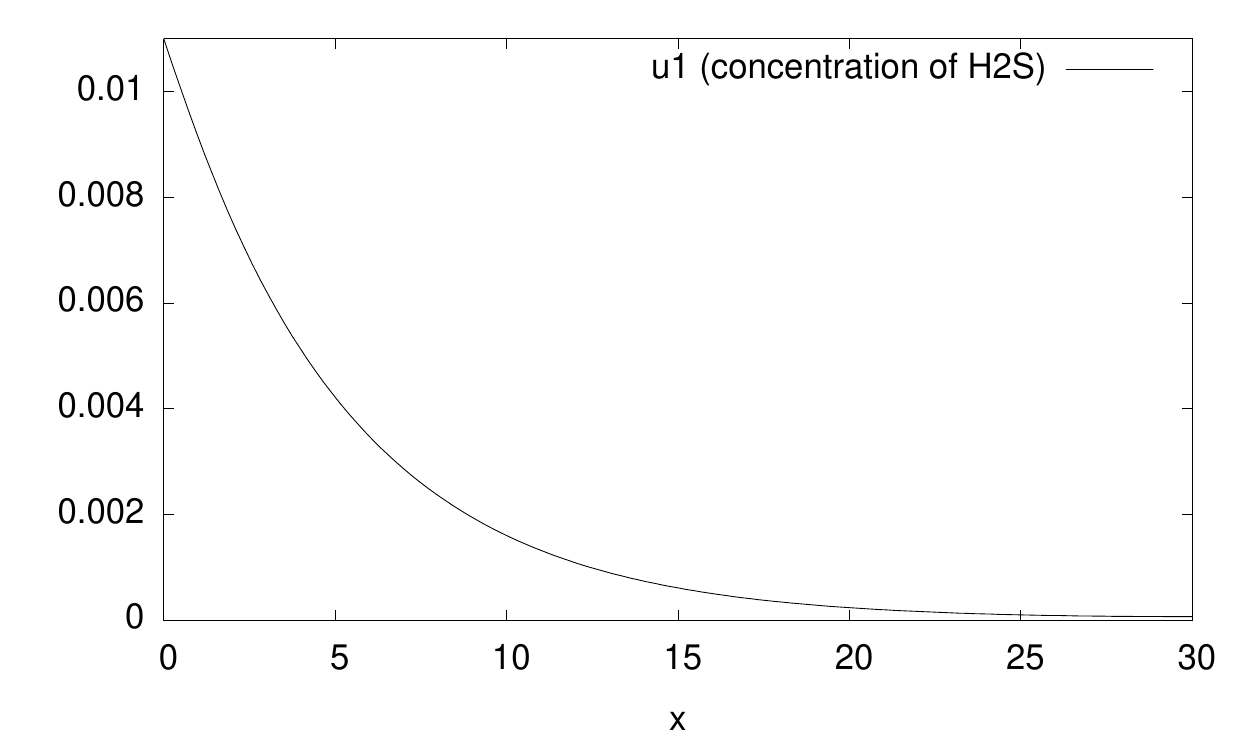}\\
    \includegraphics[width=0.48\textwidth]{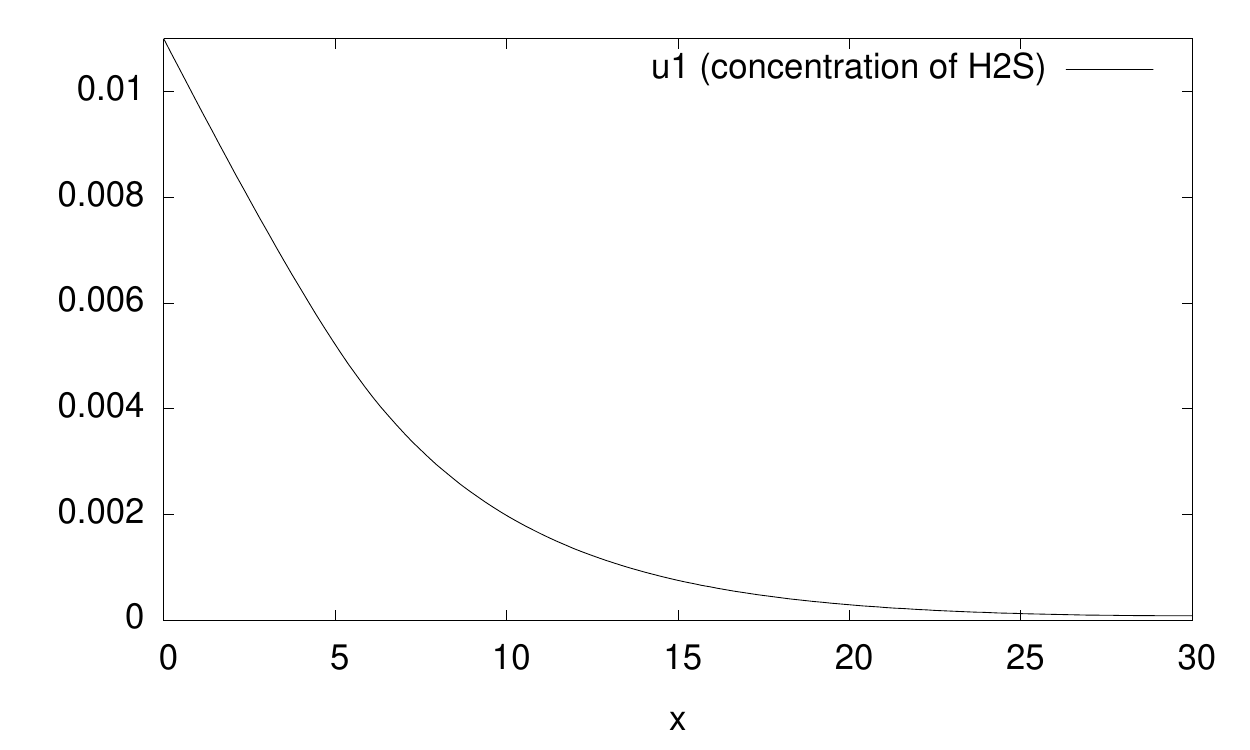}
    \includegraphics[width=0.48\textwidth]{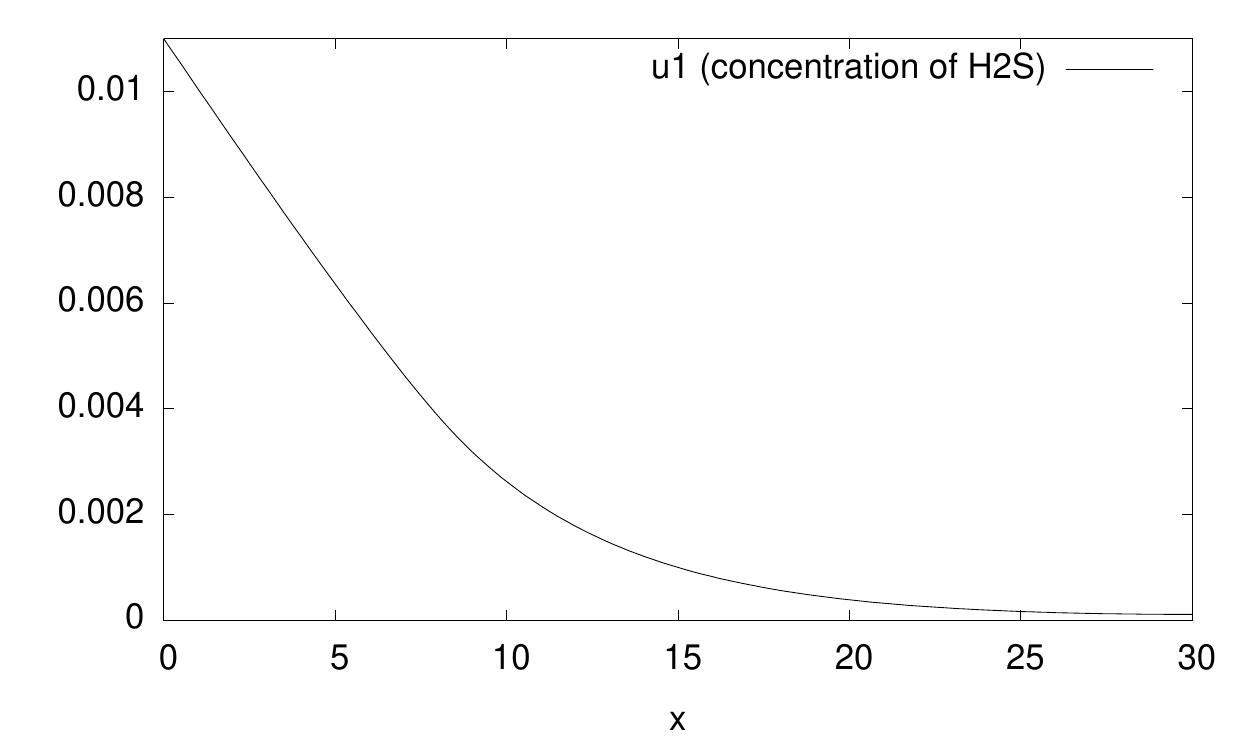}\\
    \includegraphics[width=0.48\textwidth]{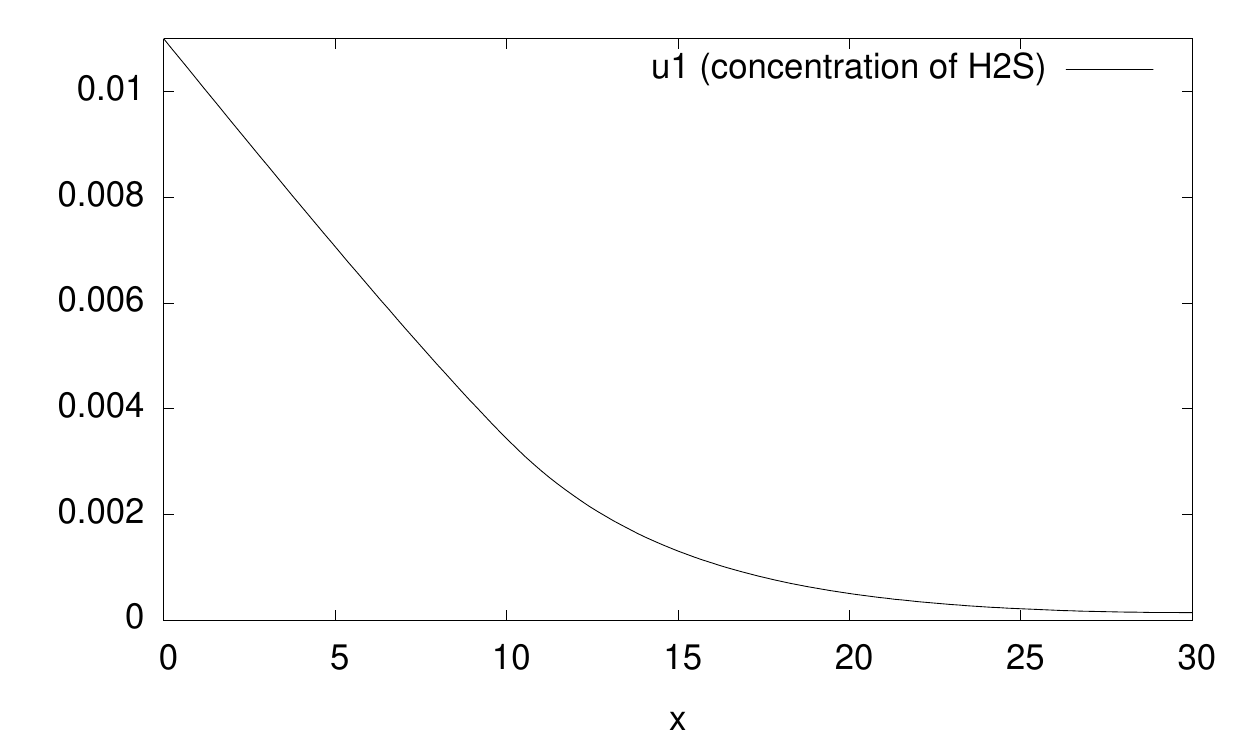}
    \includegraphics[width=0.48\textwidth]{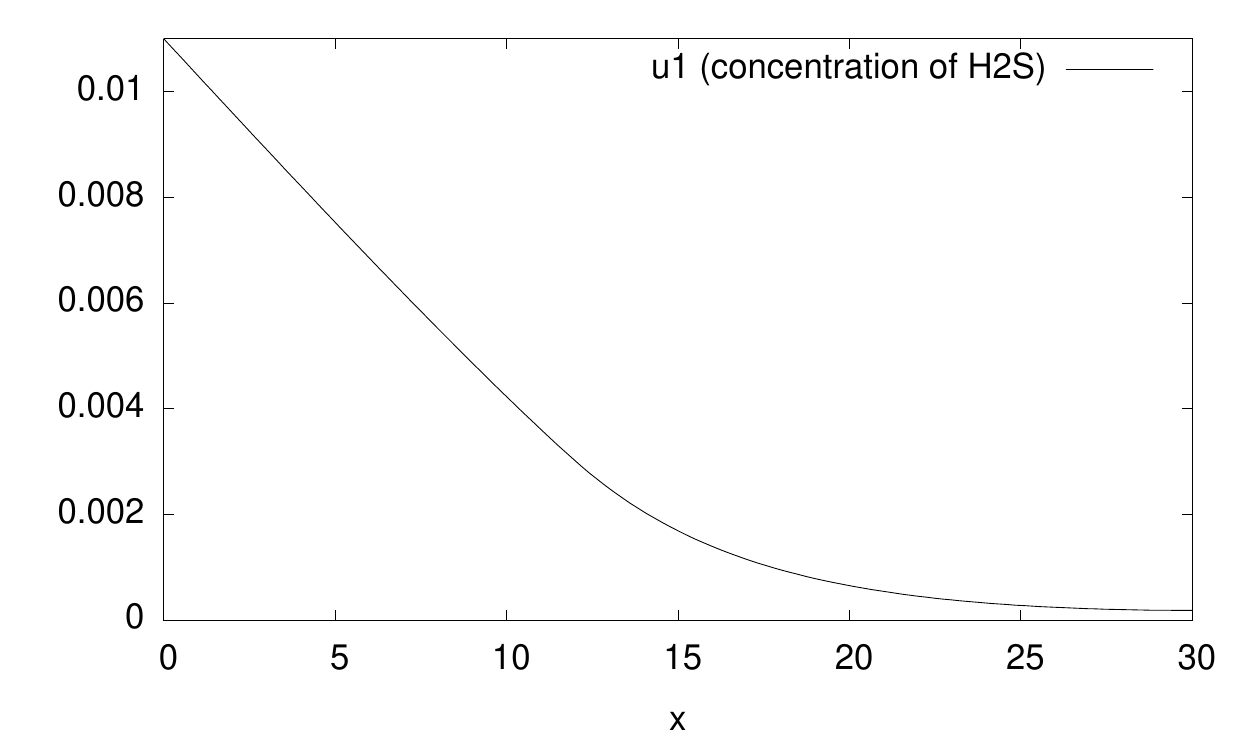}
  \end{center}
  \caption{Time evolution of $u_1$ (concentration of $\mathrm{H_2S(g)}$) shown at
  $t\in\{2000,4000,8000,12000,16000,20000\}$ in left-right and top-bottom order.}
  \label{fig:h2s}
\end{figure}

\begin{figure}[ht]
  \begin{center}
    \includegraphics[width=0.48\textwidth]{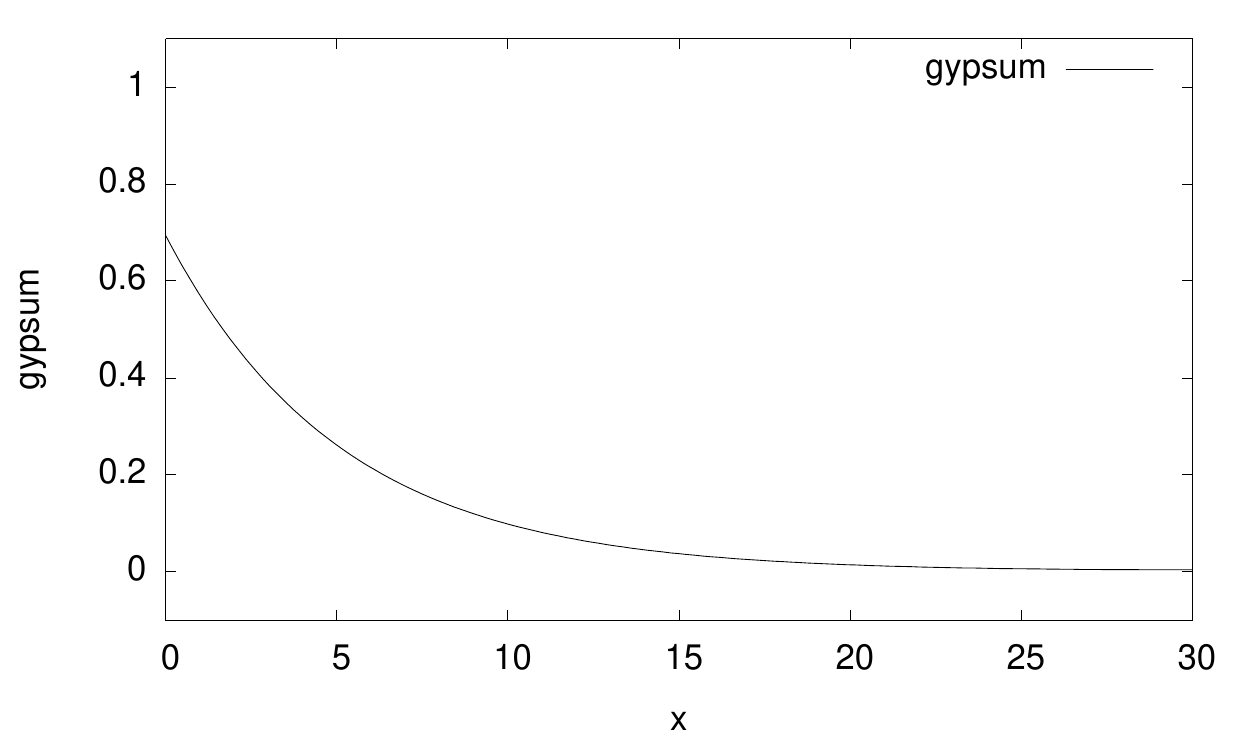}
    \includegraphics[width=0.48\textwidth]{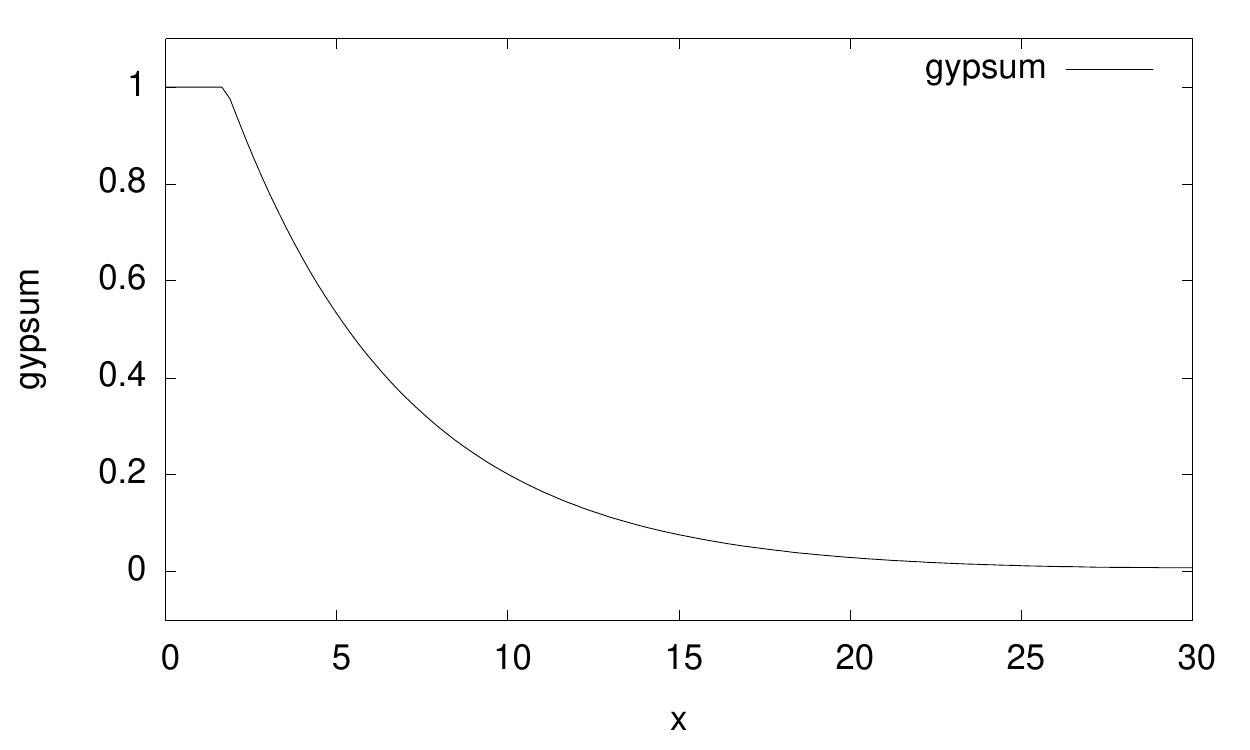}\\
    \includegraphics[width=0.48\textwidth]{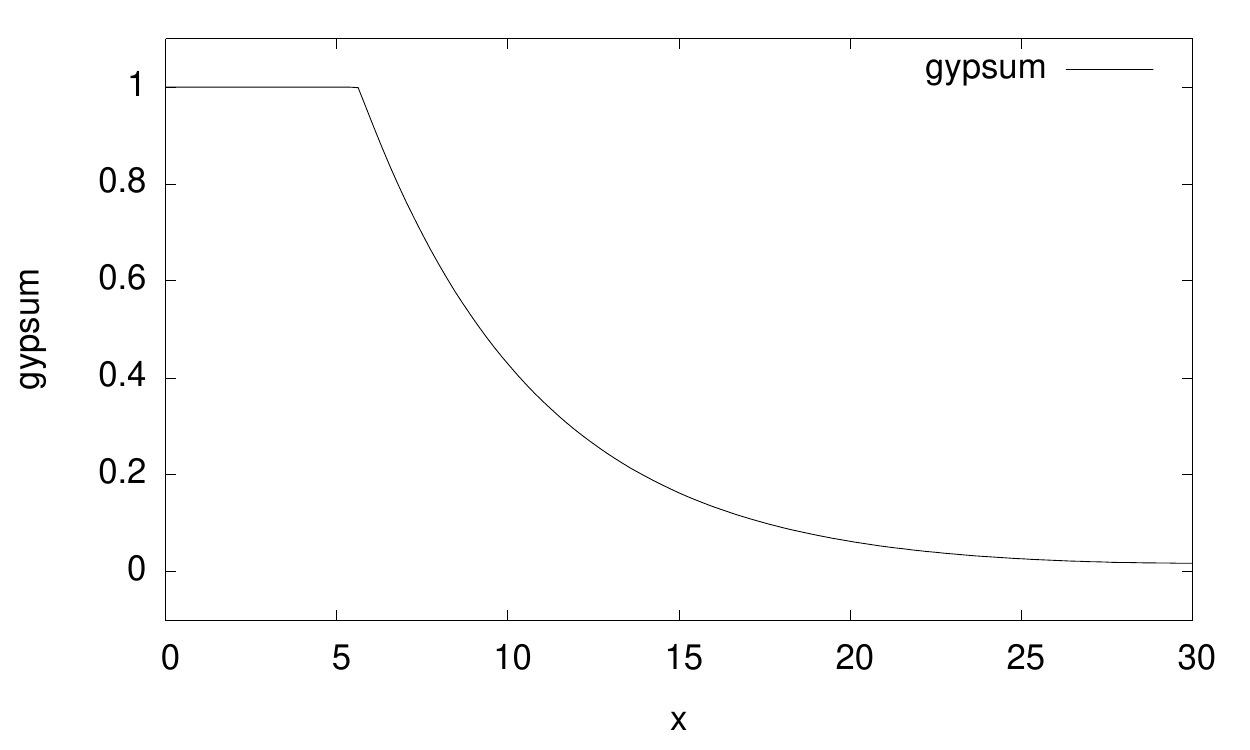}
    \includegraphics[width=0.48\textwidth]{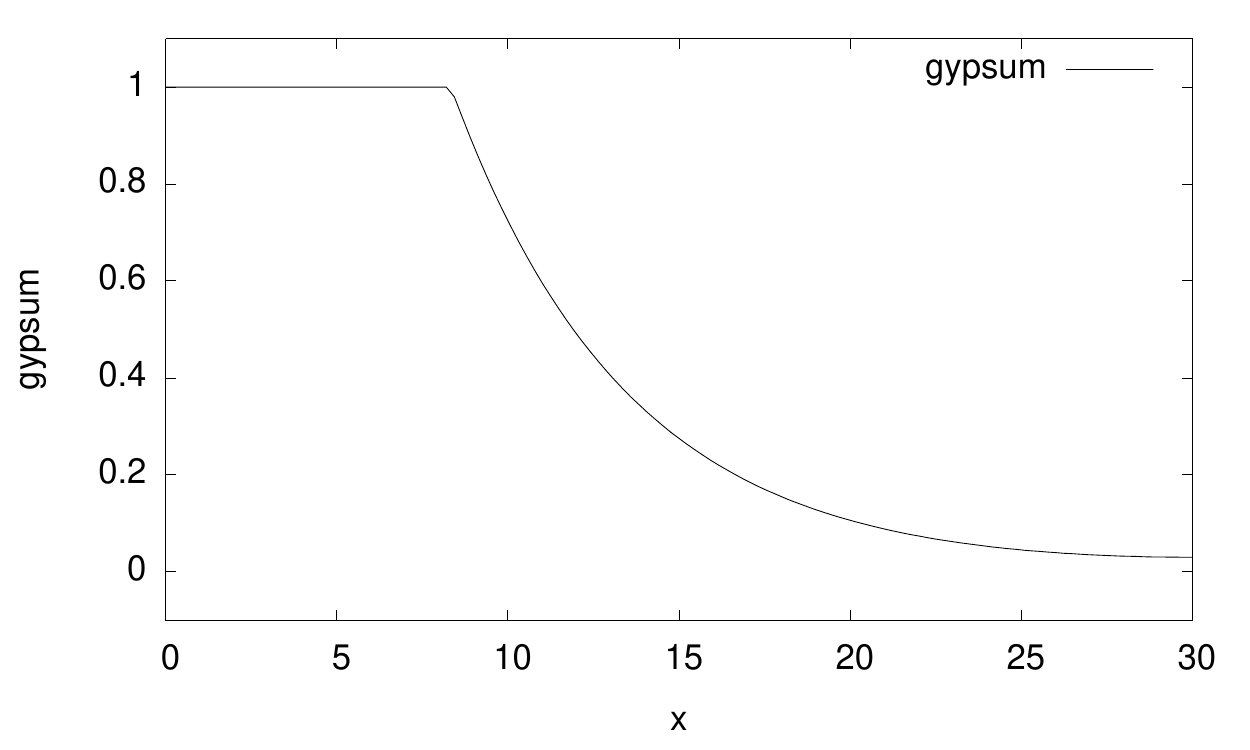}\\
    \includegraphics[width=0.48\textwidth]{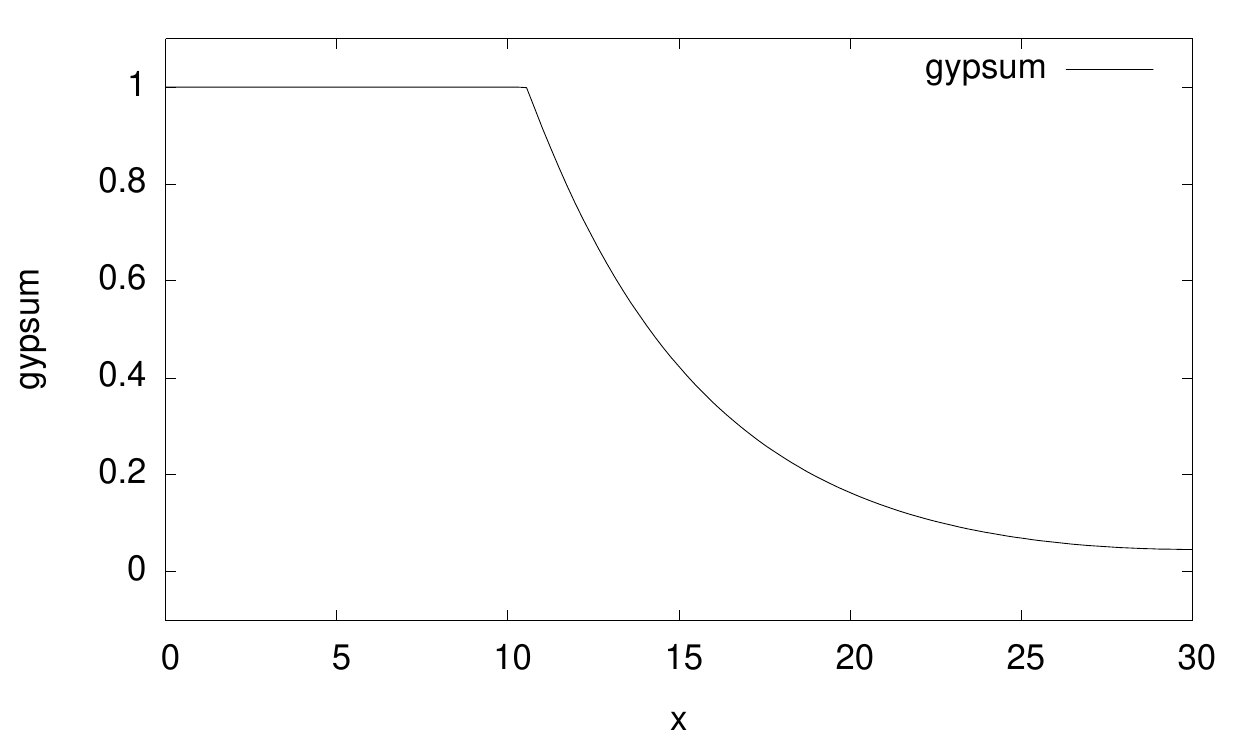}
    \includegraphics[width=0.48\textwidth]{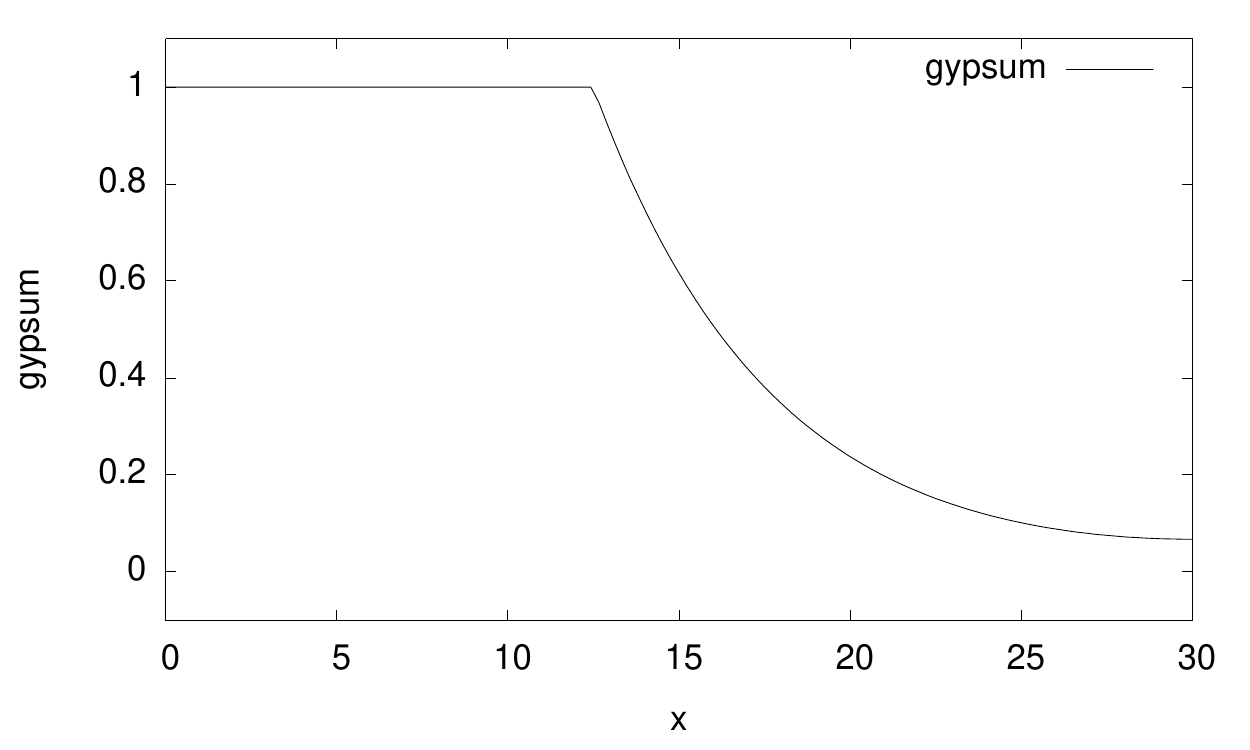}
  \end{center}
  \caption{Time evolution of $u_4$ (gypsum) shown at
  $t\in\{2000,4000,8000,12000,16000,20000\}$ in left-right and top-bottom order.}
  \label{fig:gypsum}
\end{figure}

\begin{figure}[ht]
  \begin{center}
    \includegraphics[width=0.8\textwidth]{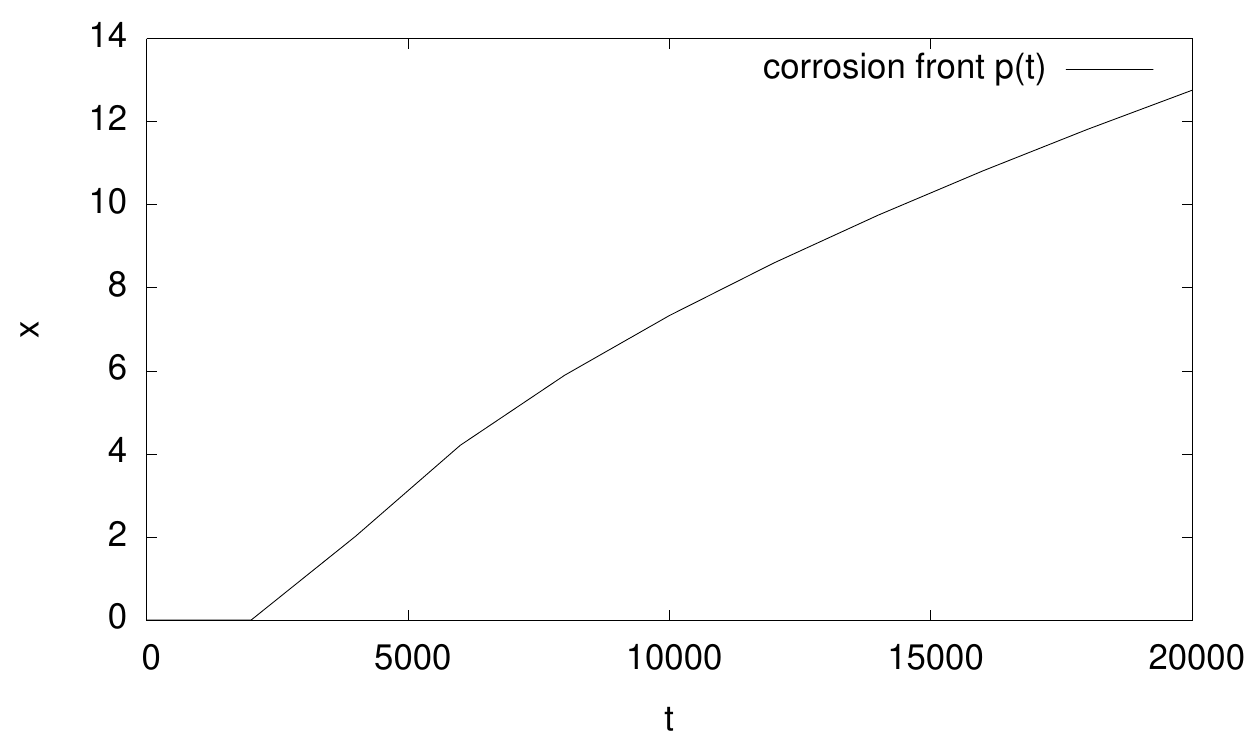}
  \end{center}
  \caption{Position of the corrosion front.}
  \label{fig:front}
\end{figure}

\begin{figure}[ht]
  \begin{center}
    \includegraphics[width=0.48\textwidth]{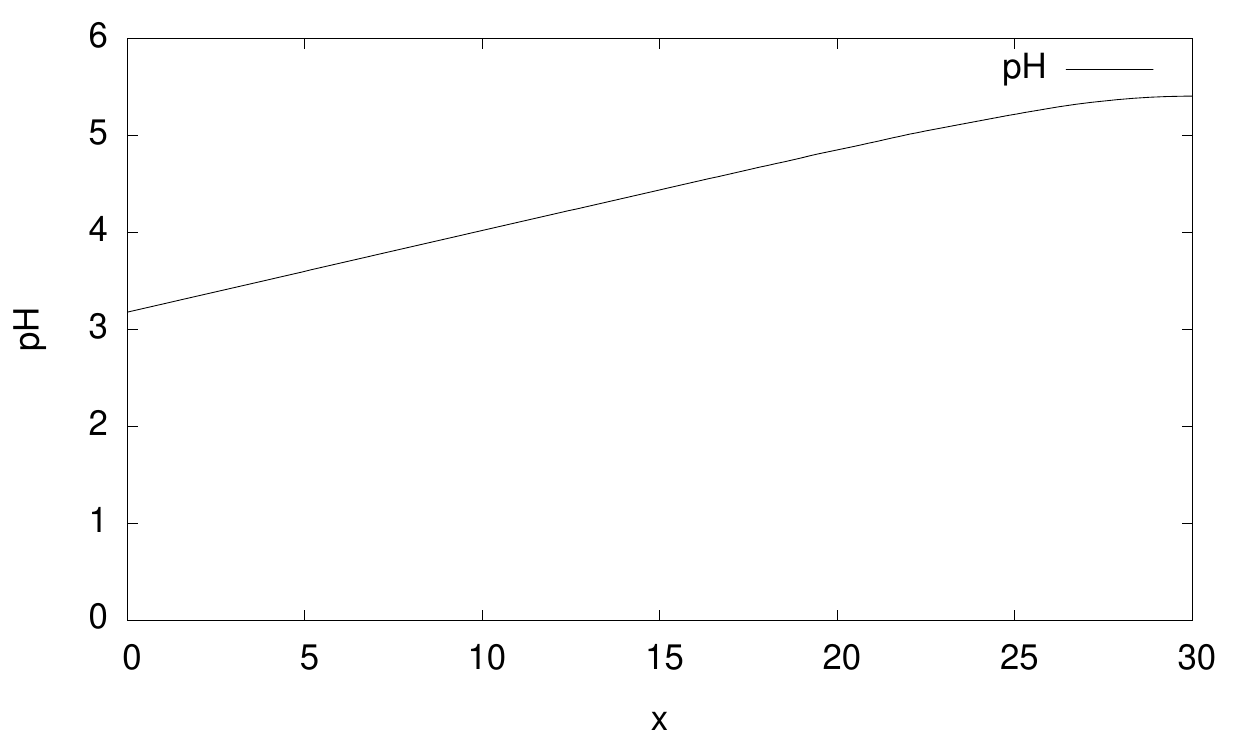}
    \includegraphics[width=0.48\textwidth]{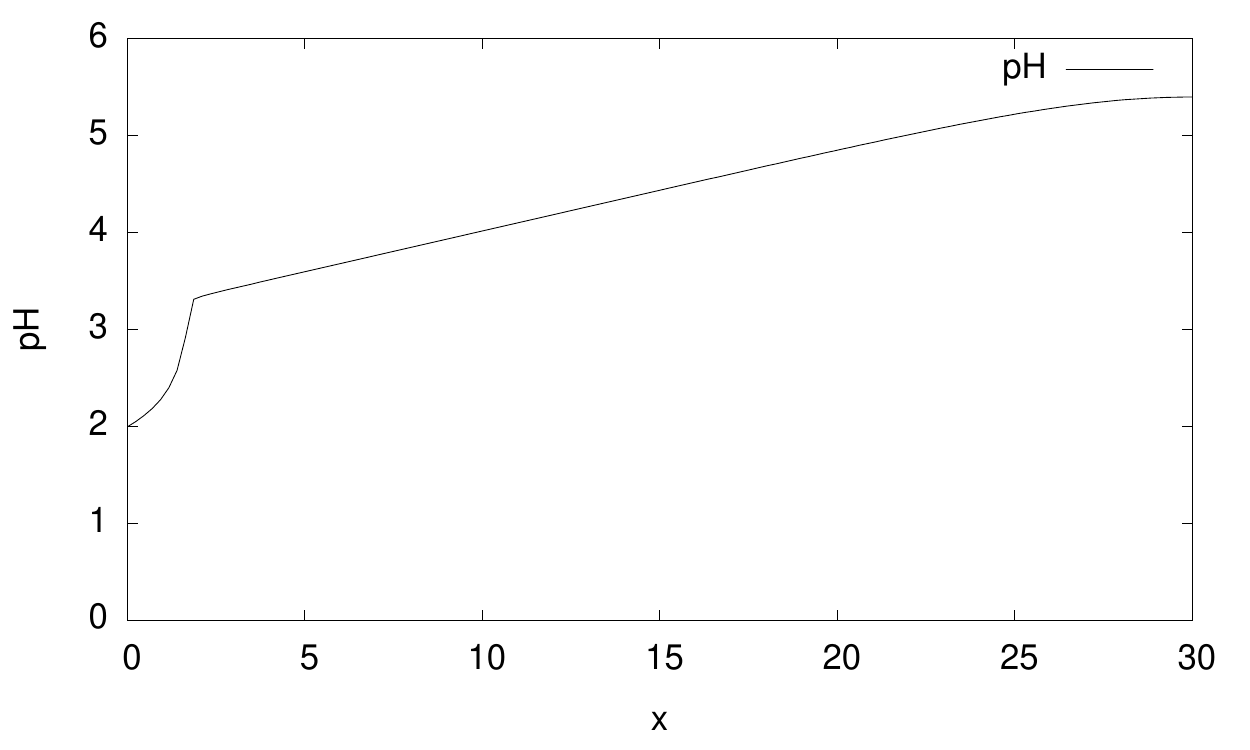}\\
    \includegraphics[width=0.48\textwidth]{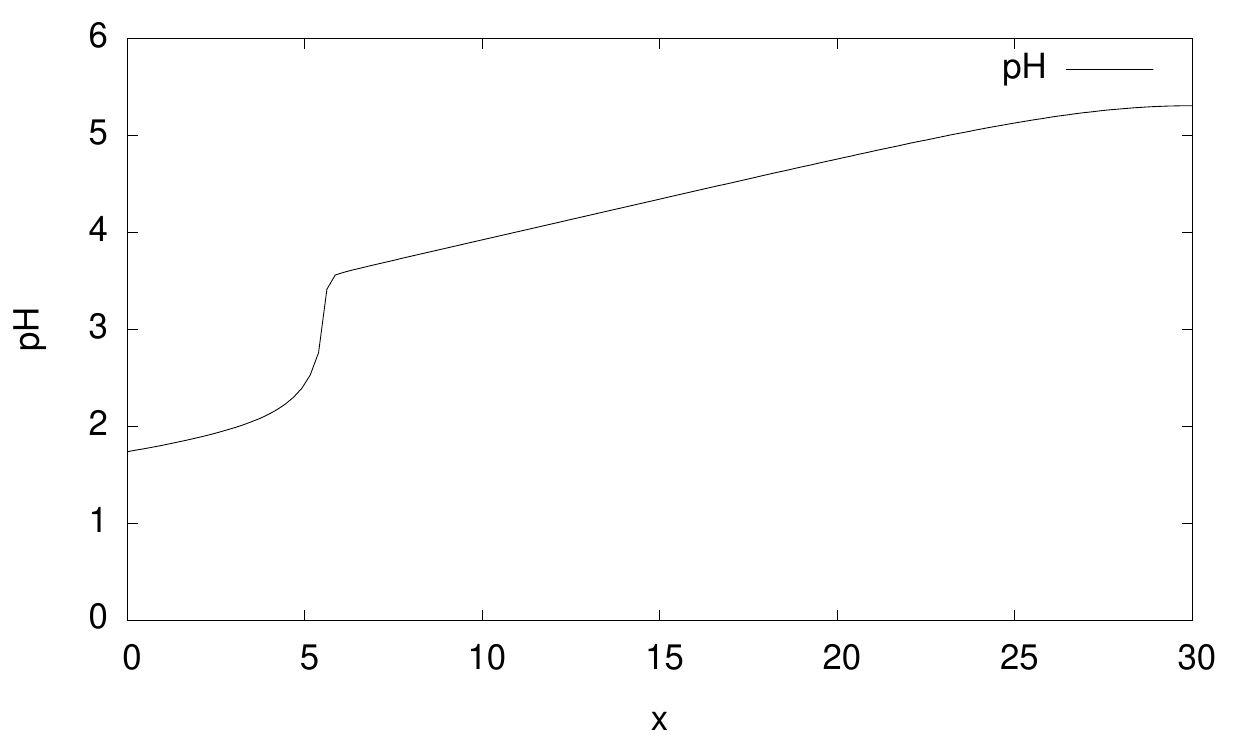}
    \includegraphics[width=0.48\textwidth]{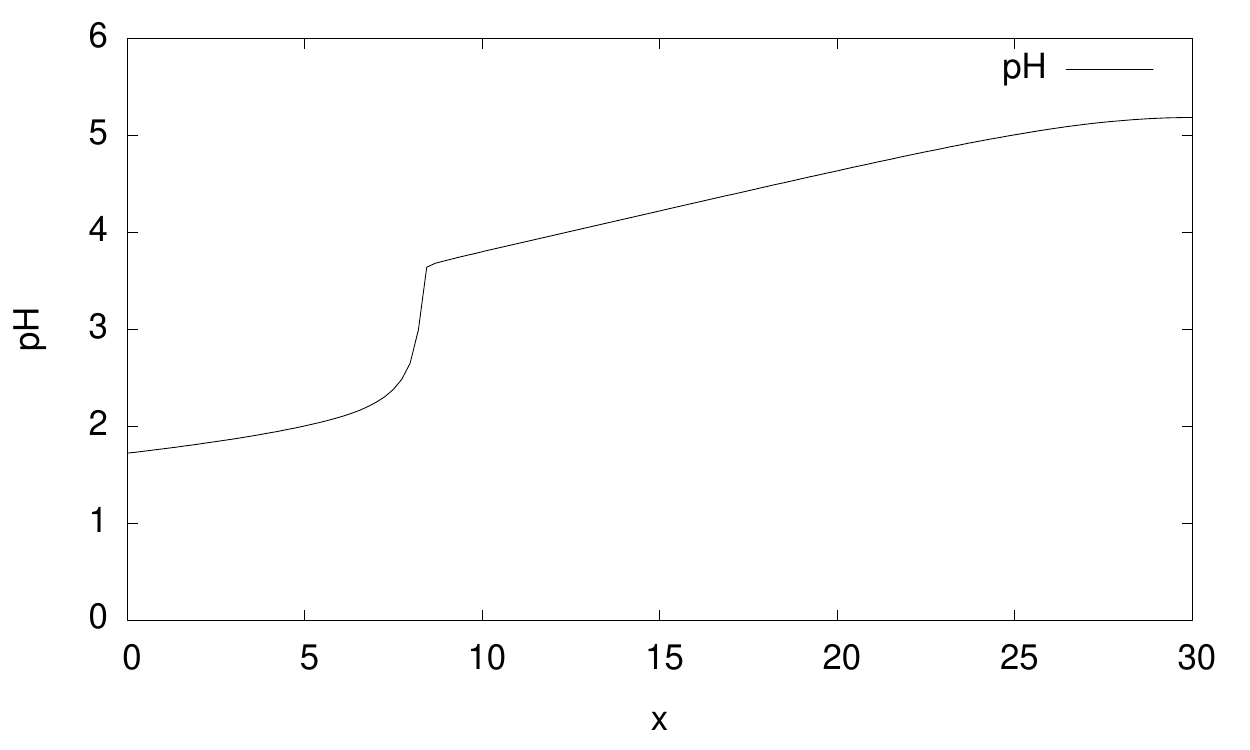}\\
    \includegraphics[width=0.48\textwidth]{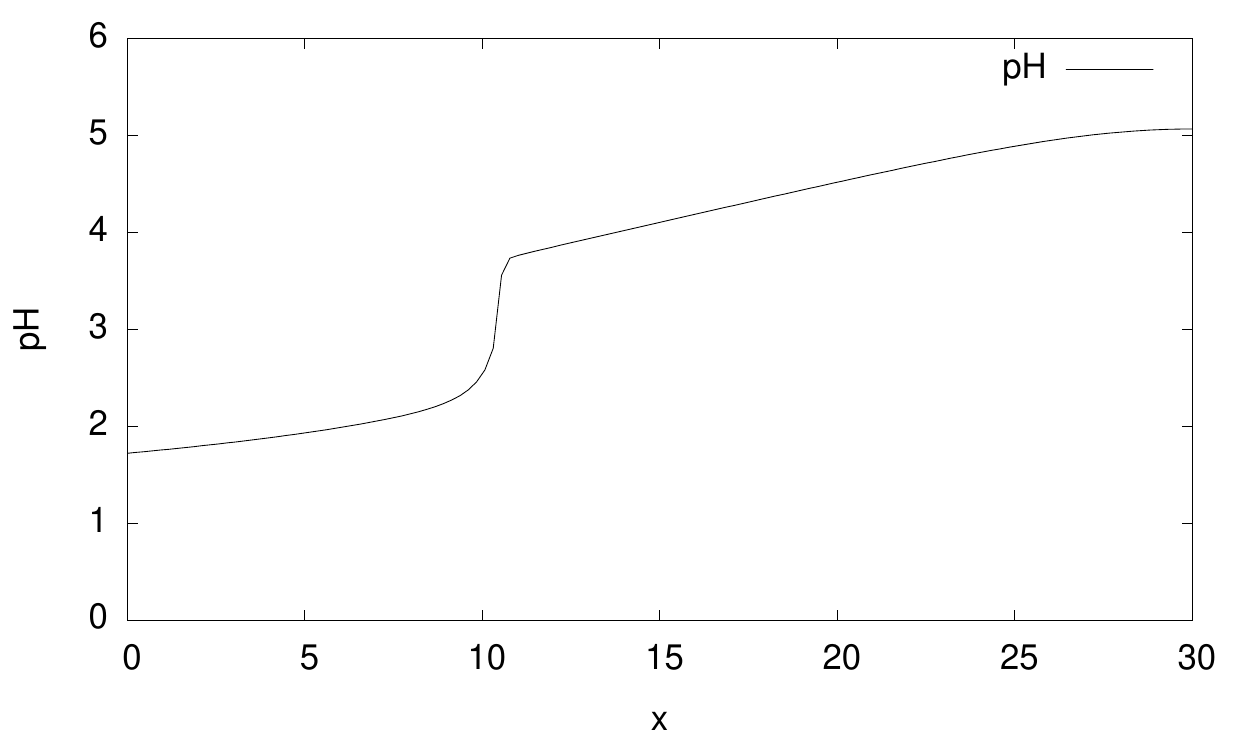}
    \includegraphics[width=0.48\textwidth]{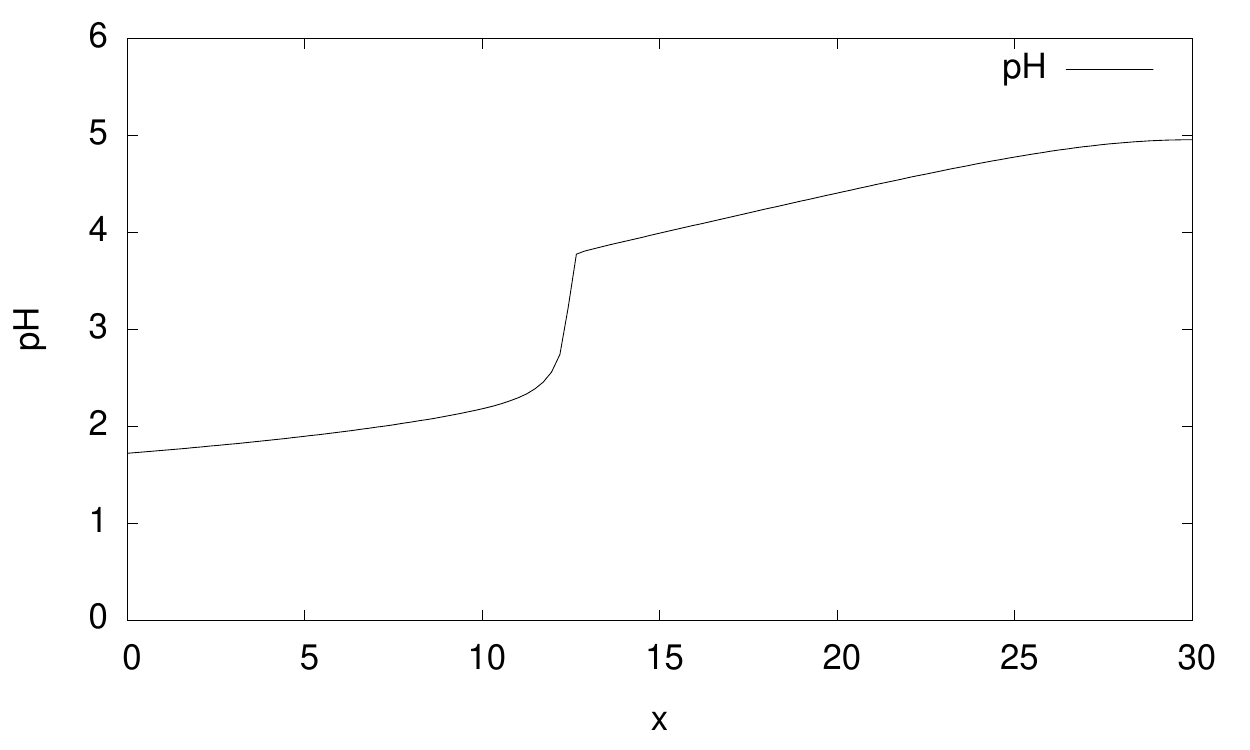}
  \end{center}
  \caption{Time evolution of macro-scale pH profiles computed from micro-scale information shown at
  $t\in\{2000,4000,8000,12000,16000,20000\}$ in left-right and top-bottom order.}
  \label{fig:ph}
\end{figure}

For the purpose of the simulation we employ a numerical scheme for a reduced
1D/2D version of the system \eqref{eq:two-scale-system}. The reduction consists
in taking $\Omega:=(0,L)$ and $Y^w=(0,\ell)$ as one-dimensional intervals, which
in effect corresponds to analysing the specimen in a perpendicular direction to
its surface away from edges and to simplifying the micro cell geometry,
respectively. The numerical scheme is based on the method of lines, where in
space we use finite difference discretization and in time we employ an implicit
higher-order time integrator for the solution of the non-linear ODE system. See
\cite{ChalupeckyEtAl2010,ChalupeckyMuntean2012} for further details of the
numerical scheme, its analysis with respect to convergence to the weak solutions
and some basic numerical experiments. For details of a computer implementation
of the numerical scheme we refer the reader to
\cite[Chapter~7]{MunteanChalupecky2011}.

In Table~\ref{tab:parameters} we summarize values of the model parameters used
in the simulations described below.

\begin{table}[h]
  \centering
  \begin{tabular}{cccccccccccccccc}
    \toprule
    $d_1$ & $d_{2,3}$ & $k_2$ & $k_3$  & $k_4$ & $\Phi_{2,3,4}$ & $B$  & $H$ & $\bar{\beta}$ & $p$ & $q$ & $u_1^D$ & $L$ & $\ell$ \\
    \midrule                                          
    0.864 & 0.00864   & 1.48  & 0.0084 & 10    & 1              & 86.4 & 0.3 & 1             & 1   & 1   & 0.011   & 30  & 1      \\
    \bottomrule
  \end{tabular}
  \caption{Parameter values used in the numerical simulation.}
  \label{tab:parameters}
\end{table}

\subsection{Free boundaries}

Figures~\ref{fig:h2s} and \ref{fig:gypsum} show the evolution of $u_1(x,t)$ and
$u_4(x,t)$ in time. The Dirichlet boundary condition $u_1(0,t)=u_1^D$ models a
constant inflow of $\mathrm{H_2S(g)}$ at $x=0$, i.e., at the surface of the
specimen. As the gas diffuses through the porous structure, it enters the water
film in the pores due to the reaction \eqref{eq:reaction-equations-3}, where it
undergoes biogenic oxidation to sulfuric acid. Consequently, its concentration
decreases with increasing depth. As the system becomes saturated and as the
sulfatation reaction \eqref{eq:sulfatation} converts available cement into
gypsum, the total concentration of $\mathrm{H_2S(g)}$ starts to increase
(Fig.~\ref{fig:h2s}).

Sulfuric acid that arises from the oxidation of $\mathrm{H_2S(aq)}$ then reacts
at $y=\ell$ with the cement paste and converts it into gypsum ($u_4$) whose
concentration profile is shown in Fig.~\ref{fig:gypsum}. Interestingly, although
the behavior of $u_1$ is as expected (i.e., purely diffusive), we notice that a
\emph{macroscopic} gypsum layer (region where $u_4$ is produced) is formed
around $t=1500$ and grows in time. The figure clearly indicates that there are
two distinct regions separated by a slowly moving intermediate layer: the left
region---the place where the gypsum production reached saturation (a threshold),
and the right region---the place of the ongoing sulfatation reaction
\eqref{eq:sulfatation} (the gypsum production has not yet reached here the
natural threshold).

We use $u_4$ to extract an approximate position of the corrosion front $p(t)$
which we define as (in our scenario, we expect $u_4$ to be decreasing)
\begin{equation*}
  p(t) := \{x\in (0,L)\ |\ u_4(x,t) = \bar{\beta} - \varepsilon\},
\end{equation*}
where $\varepsilon$ is a small parameter. Figure~\ref{fig:front} shows a graph
of $p(t)$ arising from our numerical experiment. We notice that as the corroding
front advances further into the concrete specimen, its rate of growth decreases.
This is in agreement with experimental data since the hydrogen sulfide gas
supplied from the outside environment has to be transported (by diffusion) over
ever larger distance. It is important to note that the precise position of the
separating layer is \emph{a~priori unknown} and to capture it simultaneously
with the computation of the concentration profile would require a
moving-boundary formulation similar to the one reported in \cite{mbp}.

\subsection{Drop in pH}

Emission of hydrogen sulfide from the wastewater to the air space of sewer pipe
is an important process because the problems of hydrogen sulfide in sewer pipes
are associated with gaseous hydrogen sulfide. Hydrogen sulfide is a weak acid
with a dissociation constant of 7.0 (at $20^\circ$C) and only the
non-dissociated form is emitted in the air space sewer pipe. The pH of the
wastewater is therefore of importance when evaluating the potential hydrogen
sulfide emission. After the hydrogen sulfide arrives at and diffuses into the
concrete, the oxidation of hydrogen sulfide is biological once the pH of the
solid matrix has dropped below approximately 8--9 \cite{Parker}. This represents
the tendency of hydronium ions to interact with other components of the solution,
which affects among other things the electrical potential read using a pH meter. 
The concentration of hydrogen ions is expressed as pH scale and pH is defined as
a negative decimal logarithm of the concentration of hydronium ions dissolved in
a solution.

As the model considered in this paper contains information both from micro
(pore) scale and macro-scale, we need a way of computing macro-scale values of
pH from the available micro-scale data. Such information is not readily
available in the system \eqref{eq:two-scale-system}. However, sulfuric acid is a
diprotic acid with two stages of dissociation, where the first stage occurs
fully and the dissociation in the second stage can be neglected. Therefore, the
concentration of hydronium ions is proportional to the concentration of sulfuric
acid, which is represented by $u_3$ in our model. We extract the macro-scale
concentration of sulfuric acid at each $x$ by taking a volume average of $u_3$ over
$Y^w$. Thus, we use the following expression for computing macroscopic pH:
\begin{equation}\label{eq:ph}
  \pHmac (x,t) = -\log_{10}\left(\frac{k_a}{|Y^w|} \int_{Y^w} u_3(x,y,t) \dx{y}
  \right), 
\end{equation}
where $k_a$ is the activity of hydronium ions.

The macro-scale pH profile computed using the formula \eqref{eq:ph} is shown in
Figure~\ref{fig:ph}. We can see that at the beginning of the simulation (first
graph) with increasing depth the pH also increases from acid to more basic
values as expected. Once all the available cement is consumed and converted into
gypsum (this happens for the first time at $x=0$ between the first and second
graph in Figure~\ref{fig:ph} around $t=1500$), the pH drops rapidly across the
corrosion front. This is due to the fact that behind the corrosion front the
sulfuric acid is no longer neutralised by the sulfatation reaction
\eqref{eq:sulfatation}.

Note that our pH profiles are not in the experimental range. We expect the size
of the drop will become comparable to the one seen in experiments as soon as
effects of nonlinear moisture transport and bacteria motility and chemical
activity are taken into account in the model equations. The main message that we
want to draw is that we are able to detect and compute a macroscopic pH drop,
once the needed micro-information is available.

\section*{Acknowledgments}

We thank Peter Raats (Wageningen), Hans Kuipers (Eindhoven) and Varvara
Kouznetsova (Eindhoven) for useful discussions on this topic during a multiscale
symposium that took place in Eindhoven in September 2011.

\bibliographystyle{amsplain}
\bibliography{VTJA}

\end{document}